\documentclass[12pt,oneside,english]{amsart}
\usepackage[latin1]{inputenc}
\usepackage{amssymb}

\makeatletter
\newtheorem{theorem}{Theorem}
\newtheorem{lemma}{Lemma}

\newtheorem{remark}{Remark}

\usepackage{babel}

\makeatother
\begin{document}
\baselineskip=17pt

\title[Generalized Newman Phenomena]{Generalized Newman Phenomena and Digit Conjectures on Primes}

\author{Vladimir Shevelev}
\address{Departments of Mathematics \\Ben-Gurion University of the
 Negev\\Beer-Sheva 84105, Israel. e-mail:shevelev@bgu.ac.il}


\begin{abstract}
We prove that the ratio of the Newman sum over numbers multiple of a fixed integer which is not multiple of 3 and the Newman sum over numbers multiple of a fixed integer divisible by 3 is o(1) when the upper limit of summing tends to infinity. We also discuss a connection of our results with a digit conjecture on primes.
\end{abstract}

\maketitle

\section{Introduction}

     Denote for   $x,\;\;m\in \mathbb{N}$
\begin{equation}\label{1}
S_m(x)=\sum_{0\leq n < x,n\equiv 0\pmod m}   (-1)^{s(n)}  ,
\end{equation}

where $s(n)$ is the number of 1's in the binary expansion of $n$. Sum (\ref{1}) is a \slshape Newman digit sum.\upshape

     From the fundamental paper of A.O.Gelfond \cite{4} it follows that
\begin{equation}\label{2}
S_m(x)=O(x^\lambda),\;\;\lambda=\frac{\ln 3}{\ln 4}.
\end{equation}

The case $m=3$ was studied in detail \cite{5}, \cite{2}, \cite{7}.

So, from the Coquet's theorem  \cite{2}, \cite{1} it follows that

\begin{equation}\label{3}
-\frac 1 3 +\frac{2}{\sqrt{3}}x^\lambda\leq S_3(3x) \leq \frac 1 3 +\frac{55}{3}\left(\frac{3}{65}\right)^\lambda x^\lambda
\end{equation}

with a microscopic improvement \cite{7}

\begin{equation}\label{4}
\frac{2}{\sqrt{3}}x^\lambda \leq S_3(3x)\leq \frac {55}{3}\left(\frac{3}{65}\right)^\lambda x^\lambda,\;\;\; x\geq 2,
\end{equation}

and moreover,
\newpage

\begin{equation}\label{5}
\left\lfloor 2 \left(\frac x 6\right)^\lambda\right\rfloor \leq S_3(x) \leq\left \lceil\frac{55}{3}\left(\frac{x}{65}\right)^\lambda\right\rceil.
\end{equation}

These estimates give the most exact modern limits of the so-called \slshape Newman phenomena\upshape. Note that
M.Drmota and M.Skalba \cite{3} using a close function $(S_m^{(m)}(x))$ proved that if $m$ is a multiple of 3 then for sufficiently large $x$,

\begin{equation}\label{6}
S_m(x) > 0,\;\;\;x\geq x_0(m).
\end{equation}

In this paper we study a general case for $m\geq 5$ (in the cases of $m=2$ and $m=4$ we have $|S_m(n)|\leq 1)$.

To formulate our results put for $m\geq 5$

\begin{equation}\label{7}
\lambda_m=1+\log_2b_m,
\end{equation}

\begin{equation}\label{8}
\mu_m=\frac{2b_m+1}{2b_m-1},
\end{equation}

where

\begin{equation}\label{9}
b_m^2=\begin{cases}\sin( \frac {\pi} {3}(1+\frac {3} {m})) (\sqrt{3}-\sin(\frac {\pi} {3}(1+\frac {3} {m}))) ,\;\; if \;\; m\equiv 0\pmod 3\\
\sin(\frac {\pi} {3}(1-\frac {1} {m})) (\sqrt{3}-\sin(\frac {\pi} {3}(1-\frac {1} {m}))) ,\;\; if \;\; m\equiv 1\pmod3\\
\sin(\frac {\pi} {3}(1+\frac {1} {m})) (\sqrt{3}-\sin(\frac {\pi} {3}(1+\frac {1} {m}))) ,\;\; if \;\; m\equiv 2\pmod3.\end{cases}
\end{equation}

Directly one can see that

\begin{equation}\label{10}
\frac{\sqrt{3}}{2}> b_m \geq \begin{cases} 0.861 840 88 \ldots,\;\; if \;\;(m, 3)=1,\\
0.85559967 \ldots, \;\; if \;\; (m, 3) =3,\end{cases}
\end{equation}

and thus,

\begin{equation}\label{11}
\lambda_m < \lambda
\end{equation}

and

\begin{equation}\label{12}
3.73205080\ldots < \mu_m \leq  \begin{cases} 3.76364572\ldots,\;\; if \;\;(m, 3)=1,\\
3.81215109 \ldots, \;\; if \;\; (m, 3) =3.\end{cases}
\end{equation}

Below we prove the following results.

\newpage

\begin{theorem}\label{t1}

If $(m, 3)=1$ then

\begin{equation}\label{13}
|S_m(x)|\leq \mu_m x^{\lambda_m}.
\end{equation}
\end{theorem}

\begin{theorem}\label{t2}
(Generalized Newman phenomena).  If $m>3$ is a multiple of 3 then

\begin{equation}\label{14}
\left|S_m(x)-\frac 3 m   S_3(x)\right|\leq\mu_m x^{\lambda_m}.
\end{equation}
\end{theorem}

Using Theorem 2 and  (\ref{5}) one can estimate $x_0(m)$ in (\ref{6}).   E.g., one can prove that $x_0(21)<e^{985}$.

\section{Explicit formula for $S_m(N)$}

We have

$$
S_m(N) = \sum^{N-1}_{n=0,m|n}(-1)^{s(n)}=\frac 1 m \sum^{m-1}_{t=0}\sum^{N-1}_{n=0}(-1)^{s(n)}e^{2\pi i (\frac{nt}{m})}=
$$

\begin{equation}\label{15}
=\frac 1 m \sum^{m-1}_{t=0}\sum^{N-1}_{n=0}e^{2\pi i (\frac t m n+\frac 1 2 s(n))}.
\end{equation}

Note that the interior sum has the form

\begin{equation}\label{16}
F_\alpha (N)=\sum^{N-1}_{n=0}e^{2\pi i(\alpha n + \frac 1 2 s(n))}\;\; 0\leq \alpha < 1.
\end{equation}

\begin{lemma}\label{l1}
If $N= 2^{\nu_0}+ 2^{\nu_1}+ \ldots + 2 ^{\nu_r},\;\;\nu_0 > \nu_1 > \ldots >\nu_r\geq 0$, then

\begin{equation}\label{17}
F_\alpha(N)=\sum^r_{h=0}e^{2\pi i(\alpha \sum^{h-1}_{j=0} 2^{\nu_j}+\frac h 2)}\prod^{\nu_h-1}_{k=0}(1+e^{2\pi i (\alpha 2^k +\frac 1 2)}),
\end{equation}

where as usual $\sum^{-1}_{j=0}=0,\;\;\prod^{-1}_{k=0}=1$.
\end{lemma}

\bfseries Proof. \mdseries  Let $r=0$. Then by (\ref{16})

$$
F_\alpha(N) = \sum^{N-1}_{n=0}(-1)^{s(n)}e^{2\pi i \alpha n}=1-\sum^{\nu_0-1}_{j=0}e^{2\pi i \alpha 2^j}+
$$
\begin{equation}\label{18}
+\sum_{0\leq j_1 < j_2\leq\nu_0-1}e^{2\pi i \alpha(2^{j_1}+2^{j_2})}-\ldots=\prod^{\nu_0-1}_{k=0}(1-e^{2\pi i \alpha 2^k}),
\end{equation}
\newpage

which corresponds to (\ref{17}) for $r=0$.

Assuming that (\ref{17}) is valid for every $N$ with $s(N)=r+1$ let us consider $N_1=2^{\nu_r}a + 2^{\nu_{r+1}}$ where $a$ is odd, $s(a)=r+1$ and $\nu_{r+1}<\nu_r$. Let

$$
N=2^{\nu_r}a=2^{\nu_0}+\ldots +2^{\nu_r};\;\; N_1=2^{\nu_0}+ \ldots +2^{\nu_r}+2^{\nu_{r+1}}.
$$

Notice that for $n\in[0, 2^{\nu_{r+1}})$ we have

$$
s(N+n)=s(N)+s(n).
$$

Therefore,

$$
F_\alpha(N_1) =F_\alpha(N)+\sum^{N_1-1}_{n=N}e^{2\pi i (\alpha n+\frac 1 2 s(n))}=
$$

$$
=F_\alpha(N)+\sum^{2^{\nu_{r+1}}-1}_{n=0}e^{2\pi i(\alpha n+\alpha N +\frac 1 2(s(N)+s(n)))}=
$$

$$
=F_\alpha(N)+e^{2\pi i (\alpha N+\frac 1 2 s(N))}\sum^{2^{\nu_{r+1}}-1}_{n=0}e^{2\pi i (\alpha n+\frac 1 2 s(n))}.
$$

Thus, by (\ref{17}) and (\ref{18}),

$$
F_\alpha(N_1)=\sum^r_{h=0}e^{2\pi i(\alpha \sum^{h-1}_{j=0}2^{\nu_j}+\frac h 2)}\prod^{\nu_h-1}_{k=0}(1+e^{2\pi i(\alpha 2^k+\frac 1 2)}+
$$

$$
+e^{2\pi i (\alpha\sum^r_{j=0}2^{\nu_j}+\frac{r+1}{2})}\prod^{\nu_{r+1}-1}_{k=0}\left(1+e ^{2\pi i(\alpha 2^k+\frac 1 2)} \right)=
$$

$$
=\sum^{r+1}_{h=0}e^{2\pi i (\alpha\sum^{h-1}_{j=0}2^{\nu_j}+\frac h 2)}\prod^{\nu_h-1}_{k-0}\left(1+e^{2 \pi i (\alpha 2^k +\frac 1 2)}\right).\blacksquare
$$

Formulas (\ref{15})-(\ref{17}) give an explicit expression for $S_m(N)$ as a linear combination of the products of the form

\begin{equation}\label{19}
\prod^{\nu_h-1}_{k=0}\left(1+e^{2\pi i (\alpha 2^k+\frac 1 2)}\right),\;\; \alpha=\frac t m,\;\; 0\leq t\leq m-1.
\end{equation}

\begin{remark}\label{r1}   On can extract (\ref{17}) from a very complicated general Gelfond formula \cite{4}, however, we prefer to give an independent proof.
\end{remark}
\newpage

\section{Proof of Theorem 1}

Note that in (\ref{17})

\begin{equation}\label{20}
r\leq \nu_0=\left\lfloor\frac{\ln N}{\ln 2}\right\rfloor.
\end{equation}

By Lemma 1 we have

$$
|F_\alpha(N)|\leq \sum_{\nu_h=\nu_0,\nu_1,\ldots,\nu_r}\left| \prod^{\nu_h}_{k=1}\left(1+e^{2\pi i (\alpha 2^{k-1}+\frac 1 2)} \right)\right|\leq
$$

\begin{equation}\label{21}
\leq \sum^{\nu_0}_{h=0}\left|\prod^h_{k=1}\left(1+e^{2\pi i(\alpha 2^{k-1}+\frac 1 2)} \right)\right|.
\end{equation}

Furthermore,

$$
1+e^{2\pi i( 2^{k-1}\alpha+\frac 1 2)}=2\sin(2^{k-1}\alpha\pi)(\sin(2^{k-1}\alpha\pi)-i \cos (2^{k-1}\alpha \pi))
$$

and, therefore,

\begin{equation}\label{22}
\left|1+e^{2\pi i (2^{k-1}\alpha+\frac 1 2)}  \right|\leq 2 \left| \sin(2^{k-1}\alpha \pi)\right|.
\end{equation}

According to (\ref{21}) let us estimate the product

\begin{equation}\label{23}
\prod^h_{k=1}(2|\sin (2^{k-1}\alpha\pi)|)=2^h \prod^h_{k=1}|\sin(2^{k-1}\alpha\pi)|,
\end{equation}

where by (\ref{15})

\begin{equation}\label{24}
\alpha=\frac t m,\;\; 0\leq t \leq m-1.
\end{equation}

Repeating arguments of \cite{4}, put

\begin{equation}\label{25}
\left|\sin(2^{k-1}\alpha\pi)\right|= t_k.
\end{equation}

Considering the function

\begin{equation}\label{26}
\rho(x)= 2x\sqrt{1-x^2},\;\; 0\leq x \leq 1,
\end{equation}
\newpage

we have

\begin{equation}\label{27}
t_k= 2t_{k-1}\sqrt{1-t^2_{k-1}}=\rho(t_{k-1}).
\end{equation}

Note that

\begin{equation}\label{28}
\rho'(x)= 2(\sqrt{1-x^2}-\frac{x^2}{\sqrt{1-x^2}})\leq -1
\end{equation}

for $x_0\leq x \leq 1$, where

\begin{equation}\label{29}
x_0=\frac{\sqrt{3}}{2}
\end{equation}

is the only positive root of the equation $\rho(x)=x$.

Show that either

\begin{equation}\label{30}
t_k\leq\sin\left(\frac {\pi}{m}\left\lfloor\frac m 3\right\rfloor \right)=\sin\left(\frac {\pi}{m}\left\lceil\frac {2m}{3}\right\rceil \right)=g_m<\frac{\sqrt{3}}{2}
\end{equation}

or simultaneously $t_k > g_m$ and

$$
t_k t_{k+1}\leq \max_{0\leq l \leq m-1}\left(\left|\sin \frac{l\pi}{m}\right |\left(\sqrt{3}-\left| \sin \frac {l\pi}{m}\right|\right)\right)=
$$

\begin{equation}\label{31}
=\begin{cases} \left(\sin\left(\frac \pi m \left\lfloor\frac m 3 \right\rfloor\right)\right)\left(\sqrt{3}-\sin \left (\frac \pi m\left\lfloor \frac m 3 \right\rfloor \right)\right),\;\; if \;\; m\equiv 1\pmod 3\\

\left(\sin \left( \frac \pi m\left\lceil \frac m 3 \right\rceil\right)\right)\left(\sqrt{3}-\sin\left(\frac \pi m \left\lceil \frac m 3 \right\rceil\right)\right),\;\; if \;\; m \equiv 2\pmod 3
\end{cases}=h_m< \frac 3 4.
\end{equation}

Indeed, let for a fixed values of $t\in[0, m-1]$ and $k\in[1,n]$

\begin{equation}\label{32}
t2^{k-1}\equiv l(\mod m), \;\; 0\leq l\leq m-1.
\end{equation}

Then

\begin{equation}\label{33}
t_k=\left|\sin\frac{l\pi}{m}\right|.
\end{equation}

     Now distinguish two cases:
1) $t_k\leq \frac{\sqrt{3}}{2}\;\;\;\;\; 2)t_k> \frac{\sqrt{3}}{2} $.

In case 1)

$$
 t_k= \frac{\sqrt{3}}{2}\leftrightarrows \frac{l\pi}{m}=\frac{r\pi}{3},\;\;(r,3)= 1
$$
\newpage

and since $0\leq l\leq m-1$ then

$$
m=\frac{3l}{r}, \;\;\; r=1,2.
$$

Because of  the condition $(m, 3)=1$, we have $t_k < \frac {\sqrt{3}}{2}$.

Thus, in (\ref{33})

$$
l\in\left[0,\left\lfloor\frac m 3\right\rfloor\right]\cup \left[\left\lceil\frac{2m}{3}\right\rceil,\;\; m\right]
$$

and (\ref{30}) follows.

In case 2) let $t_k>\frac{\sqrt{3}}{2}=x_0$ . For  $\varepsilon > 0$  put

\begin{equation}\label{34}
1+\varepsilon=\frac{t_k}{x_0}=\frac{2}{\sqrt{3}}\left|\sin(\pi 2^{k-1} \alpha)\right|,
\end{equation}

such that

$$
1-\varepsilon=2-\frac{2}{\sqrt{3}}\left|\sin(\pi 2^{k-1} \alpha)\right|
$$

and

\begin{equation}\label{35}
1-\varepsilon^2=\frac 4 3\left|\sin(\pi 2^{k-1} \alpha)\right|\left(\sqrt{3}-\left|\sin(\pi 2^{k-1} \alpha)\right|\right).
\end{equation}

By (\ref{27}) and (\ref{34}) we have

$$
t_{k+1}=\rho(t_k)=\rho((1+\varepsilon)x_0)=\rho(x_0) +\varepsilon x_0\rho'(c),
$$

where $c\in (x_0, (1+\varepsilon)x_0)$.

Thus, according to (\ref{28}) and taking into account that $\rho(x_0)=x_0$, we find

$$
t_{k+1}\leq x_0(1+\varepsilon)
$$

while by (\ref{34})

$$
t_k= x_0(1+\varepsilon).
$$

Now in view of (\ref{35}) and (\ref{29})

$$
t_kt_{k+1}\leq \left|\sin \pi 2^{k-1} \alpha \right|\left(\sqrt{3}-\left|\sin(\pi 2^{k-1} \alpha)\right| \right)
$$

and according to (\ref{32}),(\ref{33}) we obtain that

$$
t_kt_{k+1}\leq h_m,
$$

where $h_m$ is defined by (\ref{31}).

\newpage

Notice that from simple arguments and according to (\ref{9})

$$
g_m\leq\sqrt{h_m}=b_m.
$$

Therefore,

$$
\prod^h_{k=1}\left|\sin (\pi 2^{k-1} \alpha )\right|\leq (b_m^{\lfloor\frac h 2\rfloor})^2\leq b_m^{h-1}.
$$
 Now, by (\ref{21})- (\ref{22}), for $\alpha=\frac {t} {m},\enskip t=0, 1, ... , m-1,$ we have

$$
\left|F_{\frac {t} {m}}(N)\right|\leq \sum^{\nu_0}_{h=0}|\prod_{k=1}^{h}(1+e^{2\pi i(\alpha2^{k-1}+\frac 1 2)})|\leq
\sum^{\nu_0}_{h=0}2^{h}\prod_{k=1}^{h}|\sin(2^{k-1}\alpha\pi)|\leq$$
$$1+2\sum^{\nu_0}_{h=1}(2b_m)^{h-1}\leq1+2\frac {(2b_m)^{\nu_0}} {2b_m-1}.
$$

Note that, according to (\ref{7}) and (\ref{20})
$$ (2b_m)^{\nu_0}=2^{\lambda_m\nu_0}\leq 2^{\lambda_m\log_2N}=N^{\lambda_m}.$$
Thus,
$$|F_{\frac {t} {m}}(N)|\leq1+\frac {2} {2b_m-1} N^{\lambda_m}\leq\frac {2} {2b_m-1-\gamma_m} N^{\lambda_m},$$
where $\gamma_m$ is defined by the equality
$$\frac {1} {2b_m-1-\gamma_m}-\frac {1} {2b_m-1}=\frac 1 2.$$
Hence, we find
$$\gamma_m=\frac {(2b_m-1)^2} {2b_m+1}$$
and, consequently, by (\ref{8}),
$$|F_{\frac {t} {m}}(N)|\leq \frac {2b_m+1} {2b_m-1}N^{\lambda_m}=\mu_m N^{\lambda_m}.$$
Thus, the theorem follows from (\ref{15}).$\blacksquare$
\section{Proof of Theorem 2.}

Select in (\ref{15}) the summands which correspond to $t=0,\;\frac m 3,\;\frac{2m}{3}$.

We have

$$
m S_m(N)=\sum^{N-1}_{n=0}\left(e^{\pi i s(n)}+e^{2\pi i(\frac n 3 +\frac 1 2 s(n))}+ e^{2\pi i(\frac{2n}{3} + \frac 1 2 s(n))} \right)+
$$

\begin{equation}\label{36}
+\sum^{m-1}_{t=1,t\neq \frac m 3,\frac{2m}{3}}\;\;\sum^{N-1}_{n=0}e^{2\pi i(\frac t m n + \frac 1 2 s(n))}.
\end{equation}
\newpage

Since the chosen summands do not depend on $m$ and for $m=3$ the latter sum is empty then we find

\begin{equation}\label{37}
m S_m(N)= 3S_3(N)+\sum^{m-1}_{t=1,t\neq \frac m 3,\frac{2m}{3}}\;\;\sum^{N-1}_{n=0}e^{2\pi i(\frac t m n + \frac 1 2 s(n))}.
\end{equation}

Further, the last double sum is estimated by the same way as in Section 3 such that

\begin{equation}\label{38}
\left |S_m(N) - \frac 3 m S_3(N)\right |  \leq \mu_m N^{\lambda_m} \blacksquare.
\end{equation}

\begin{remark}\label{r}
Notice that from elementary arguments it follows that if $m\geq 5$ is a multiple of $3$ then

$$
\left(\sin \frac \pi m \left\lfloor\frac{m-1}{3}\right\rfloor\right)\left(\sqrt{3}-\sin \frac \pi m \left\lfloor\frac{m-1}{3}\right\rfloor\right)\leq
$$
$$
\leq\left(\sin \frac \pi m \left\lceil\frac{m+1}{3}\right\rceil\right)\left(\sqrt{3}-\sin \frac \pi m \left\lceil\frac{m+1}{3}\right\rceil\right).
$$
The latter expression is the value $of$  $b_m^2$ in this case (see (\ref{9})).
\end{remark}
\bfseries Example.\mdseries \enskip Let us find some $x_0$ such that $ S_{21}(x)>0$ for $x\geq x_0.$\newline
Supposing that $x$ is multiple of 3 and using (\ref{4}) we obtain that
$$S_3(x)\geq\frac {2} {3^{\lambda+\frac 1 2}} x^{\lambda}.$$
Therefore, putting  $ m=21$ in  Theorem 2,  we have
$$ S_{21}(x)\geq\frac 1 7 S_3(x)-\mu_{21}x^{\lambda_{21}}\geq \frac {2} {7\cdot3^{\lambda+\frac 1 2}} x^{\lambda}-\mu_{21}x^{\lambda_{21}}.$$
 Now, calculating $\lambda$ and $\lambda_m$ by (\ref {2}) and (\ref {8}), we find a required $x_0:$
$$x_0=(3.5 \cdot 3^{\lambda+\frac 1 2}\mu_{21})^{\frac {1} {\lambda-\lambda_{21}}}=e^{984.839...}.$$
\bfseries Corollary.\mdseries \slshape \enskip For $m$ which is not a multiple of 3, denote $U_m(x)$ the set of the positive integers not exceeding $x$ which are multiples of $m$ and not multiples of 3. Then \upshape

$$
\sum_{n\in U_m(x)}(-1)^{s(n)}=-\frac 1 m S_3(x)+O(x^{\lambda_m}).
$$

In particular, for sufficiently large $x$ we have
\newpage
$$
\sum_{n\in U_m(x)}(-1)^{s(n)}< 0.
$$

\bfseries Proof.\mdseries \enskip  Since

$$
|U_m(x)|=S_m(x)-S_{3m}(x)
$$

then the corollary immediately follows from Theorems 1, 2.

\section{On Newman sum over primes}

In \cite{6} we put the following binary digit conjectures on primes.

\bfseries Conjecture 1.\mdseries \slshape For all $n\in \mathbb{N},\;\;n\neq 5,6$

$$
\sum_{p\leq n}(-1)^{s(p)}\leq 0,
$$

where the summing is over all primes not exceeding $n$\upshape.

Moreover, by the observations, $\sum_{p\leq n}(-1)^{s(p)}< 0$ beginning with $n=31$.

\bfseries Conjecture 2.\mdseries

$$
\lim_{n\rightarrow\infty}\frac{\ln\left(-\sum_{p\leq n}(-1)^{s(p)}\right)}{\ln n}=\frac{\ln 3}{\ln 4}.
$$

\;\;\;\;\;\;

A heuristic proof of Conjecture 2 was given in \cite{8}. For a prime $p$, denote $V_p(x)$ the set of positive integers not exceeding $x$ for which $p$ is the least prime divisor. Show that the correctness of Conjectures 1 (for $n\geq n_0$) follows from the following very plausible statement, especially in view of the above estimates.

\bfseries Conjecture 3.\mdseries \slshape \enskip For sufficiently large $n$ we have \upshape

\begin{equation}\label{39}
\left|\sum_{5\leq p \leq\sqrt{n}}\;\sum_{j\in V_p(n),j>p}(-1)^{s(j)}\right|<\sum_{j\in V_3(n)}(-1)^{s(j)}=S_3(n)-S_6(n).
\end{equation}

Indeed, in the "worst case" (really is not satisfied) in which for all $n\geq p^2$

\begin{equation}\label{40}
\sum_{j\in V_p(n),j>p}(-1)^{s(j)}< 0,\;\;\;p\geq 5.
\end{equation}

we have a decreasing but \slshape positive \upshape sequence of sums
\newpage
$$
\sum_{j\in V_3(n),j> 3} (-1)^{s(j)},\;\; \sum_{j\in V_3(n),j> 3} (-1)^{s(j)}+\sum_{j\in V_5(n),j> 5} (-1)^{s(j)},
$$
$$
\ldots, \sum_{j\in V_3(n),j> 3} (-1)^{s(j)}+\sum_{5\leq p<\sqrt{n}}\;\sum_{j\in V_p(n),j>p} (-1)^{s(j)}>0.
$$

\;\;\;

 Hence, the "balance condition" for odd numbers \cite{8}

\begin{equation}\label{41}
\left|\sum_{j\leq n, j \;is\; odd}(-1)^{s(j)} \right|\leq 1
\end{equation}

must be ensured permanently by the excess of the odious primes. This explains Conjecture 1.

It is very interesting that for some primes $p$ most likely indeed (\ref{40}) is satisfied for all $n\geq p^2$. Such primes we call "resonance primes". Our numerous observations show that all resonance primes not exceeding $1000$ are:

$$
11,19,41,67,107,173,179,181,307,313,421,431,433,587,
$$
$$
601,631,641,647,727,787.
$$

In conclusion, note that for $p\geq 3$ we have

\begin{equation}\label{42}
\lim_{n\rightarrow\infty}\frac{|V_p(n)|}{n}= \frac 1 p \prod_{2\leq q < p}\left(1-\frac 1 q\right)
\end{equation}

such that

\begin{equation}\label{43}
\lim_{n\rightarrow\infty}\left(\sum_{p\geq 3}\frac{|V_p(n)|}{n}\right)= \frac 1 2 .
\end{equation}

Thus, using Theorems 1, 2 in the form

\begin{equation}\label{44}
S_m(n)=\begin{cases} o(S_3(n)),\;\;(m,3)=1 \\\frac 3 m S_3(n)(1+o(1)),\;\; 3| m \end{cases}
\end{equation}

and inclusion-exclusion for $p\geq 5$, we find

$$\sum_{j\in V_p(n)}(-1)^{\sigma(j)}=-\frac {3} {3p} \prod_{2\leq q < p,\enskip q\neq3}(1-\frac {1} {q})S_3(n)(1+o(1))=$$

\begin{equation}\label{45}
-\frac {3} {2p}\prod_{2\leq q < p}(1-\frac {1} {q})S_3(n)(1+o(1)).
\end{equation}
\newpage
Now in view of (\ref{5}) we obtain the following absolute result as an approximation of Conjectures 1, 2.

\begin{theorem}\label{t3}

For every prime number $p\geq5$ and sufficiently large $n\geq n_p$ we have

$$
\sum_{j\in V_p(n)}(-1)^{s(j)}< 0
$$

and, moreover,

$$
\lim_{n\rightarrow\infty}\frac{\ln(-\sum_{j\in V_p(n)}(-1)^{s(j)})}{\ln n}=\frac{\ln 3}{\ln 4}.
$$
\end{theorem}

\;\;\;\;\;\;\;

\end{document}